\documentstyle[sprocl]{article}
\input{psfig}
\def\be{\begin{equation}}
\def\ee{\end{equation}}
\def\bea{\begin{eqnarray}}
\def\eea{\end{eqnarray}}

\begin{document}
\title{NON--MINIMAL $q$--DEFORMATIONS AND ORTHOGONAL
SYMMETRIES: ${\cal U}_{\displaystyle{q}}$(SO(5)) EXAMPLE}
\author{B. Abdesselam$^{\dagger}$,
D. Arnaudon $^{\ddagger}$ and
A. Chakrabarti$^{\dagger,}$
\footnote{ \it Talk presented at the Nankai workshop, Tianjin, 1995 by
A. Chakrabarti.}}
\address{$^{\dagger}$Centre de Physique Th{\'e}orique,
Ecole Polytechnique,
91128 Palaiseau Cedex, France.
Laboratoire Propre du CNRS UPR A.0014}
\address{$^{\ddagger}$ENSLAPP \footnote{\it URA 14-36 du CNRS, associ\'ee \`a
l'E.N.S. de Lyon, et au L.A.P.P. d'Annecy-le-Vieux.},
Chemin de Bellevue BP 110,
74941 Annecy-le-Vieux Cedex, France.}


\maketitle
\abstracts{Non--minimal $q$-deformations are defined. Their role in the
explicit construction of the matrix elements of the generators of
${\cal U}_{q}(SO(5))$ on suitably parametrized bases are
exhibited. The implications are discussed.}

Symmetry is one theme of this workshop. Suppose one $q$-deforms some classical
symmetry (unitary, orthogonal,$\cdots$) intending to explore the possibilities
of applications of the symmetry thus generalized. Given such a goal, one
should go further than the formally $q$-deformed Hopf algebra. One should
construct explictly the representations irreducible ones to start with but
also, so far as possible, non-decomposable ones for $q$ a root of unity. By
explicit construction I mean a complete set of suitably parametrized basis
states spanning the space of the representation in question and the matrix
elements of the generators acting on these state vectors. The invariant
parameters and the variable indices labelling the states will each (some
very directly while others less so) have their specific significance in the
description of the phenomenon studied. the matrix elements will measure the
response of the states to the constraints of the symmetry (the action of the
generators). They will also yield the values of the crucial invariants. Unless
all these elements are obtained a central problem remains unsolved. One is not
fully equipped to explore possible applications.

As one proceeds with this program one encounters, among others, the following
fact. The $q$-deformed unitary algebras are relatively docile while the
corresponding orthogonal ones are suprisingly refractory. To give this
statement more precise content let me introduce at this point some definitions
and terminology.

Let us start with a classical quantity $x$, typically a factor in some
classical matrix element.
\vskip 0.25cm
{\it Minimal $q$-deformation:}
\begin{equation}
\begin{array}{lll}
q=1 & &\;\;\;\;\;\;\;\;\;\;\;\;q\not = 1 \\
x & \rightarrow &\;\;\;\;\;\;\;\;\;\;\;\;[x]_{p}\equiv {q^{px}-q^{-px}
\over q^{p}-q^{-p}} \\
\end{array}
\end{equation}
One may further refine this by defining the deformation to be strictly
minimal only for $p=1$ (when the subscript $p=1$ will be omitted) and to be
pseudo-minimal for $p\not = 1$.
\vskip 0.25cm
{\it Non-minimal $q$-deformation:}

An unlimited number of more complicated deformations (retaining the symmetry
$q\rightleftharpoons q^{-1}$ and the same classical limit $x$) is evidently
possible.
\vskip 0.25cm
{\it Example. 1.}
\begin{equation}
\begin{array}{l}
x \rightarrow [x_{1}]_{p_{1}}-[x_{2}]_{p_{2}},\;\;\;\;\;\;\;\;\;
(x_{1}-x_{2}=x)
\end{array}
\end{equation}
\vskip 0.25cm
{\it Example. 2.}
\begin{equation}
\begin{array}{l}
x \rightarrow [x]_{p}{[y]_{p_{1}} \over [y]_{p_{2}}}
\end{array}
\end{equation}
Apart from simple $q$-factors ($x\rightarrow [x]_{p}\;q^{f(x)}$, $f(x)$ being
some non-singular function of $x$ and possibly other parameters) we have, as
yet, encountered only the types $(2)$ and $(3)$ (but possibly with more than
one $y$-factors in $(3)$). Further study may permit the classification of all
the relevant ones.

Let us now go back to our initial statement. The well-known Gelfand--Zetlin
matrix elements $[1]$ for irreps. of $SU(N)$ are square roots of ratios of
products of integer factors. {\it A minimal $q$-deformation of each factor}
\begin{equation}
\begin{array}{l}
\left ( {x_{1}\;x_{2}\;x_{3}\cdots \over y_{1}\;y_{2}\;y_{3}\cdots }
\right )^{1/2} \rightarrow \left ( {[x_{1}]\;[x_{2}]\;[x_{3}]\cdots \over
[y_{1}]\;[y_{2}]\;[y_{3}]\cdots } \right )^{1/2}
\end{array}
\end{equation}
gives the corresponding element for ${\cal U}_{q}(SU(N))$ for generic $q$.
For $q$ a root of unity periodic representations are obtained by introducing
suitable {\it fractional parts} for each $x$ and $y$ $[2]$. Relatively simple
modifications yield other classes of representations $[3,\;4]$. One can of
course introduce a unitary transformation after deforming as in $(4)$ to
obtain complicated matrix elements. But all representations, at least for real
positive $q$, can be obtained as in $(4)$. Transformations can only introduce
spurious non-minimalities masking the basic simplicity.

In $q$-deforming $\;$orthogonal algebras the $\;$well-known $\;$prescriptions
for
${\cal U}_{q}(SU(2))$ suffice for $SO(3)(\approx SU(2))$ and
$SO(4)\approx (SU(2) \otimes SU(2))$. But even for $SO(4)$ problems (and
non-minimalities) arise if one tries to $q$-deform directly the canonical
Gelfand--Zetlin matrix elements $[5]$. The first intrinsically non-trivial
case is ${\cal U}_{q}(SO(5))$ which we discuss here showing
exactly where and how non-minimalities enter and analysing their implications.

${\cal U}_{q}(SO(5))$: Corresponding to the two unequal roots
one has two $q$-deformed Chevalley triplets. The standard Drinfeld-Jimbo Hopf
algebra ( omitting the coproducts, counits and antipodes ) becomes in our
conventions $[6,\;7]$,
\begin{equation}
\begin{array}{ll}
q^{\pm h_{1}} e_{1} = q^{\pm 1} e_{1} q^{\pm h_{1}}, &
q^{\pm h_{1}} f_{1} = q^{\mp 1} f_{1} q^{\pm h_{1}}, \\

q^{\pm 2 h_{2}} e_{1} = q^{\mp 1} e_{1} q^{\pm 2h_{2}}, &
q^{\pm 2 h_{2}} f_{1} = q^{\pm 1} f_{1} q^{\pm 2h_{2}}, \\

q^{\pm h_{1}} e_{2} = q^{\mp 1} e_{2} q^{\pm h_{1}}, &
q^{\pm h_{1}} f_{2} = q^{\pm 1} f_{2} q^{\pm h_{1}}, \\

q^{\pm h_{2}} e_{2} = q^{\pm 1} e_{2} q^{\pm h_{2}}, &
q^{\pm h_{2}} f_{2} = q^{\mp 1} f_{2} q^{\pm h_{2}}, \\

[e_{1} , f_{2}] = 0, & [e_{2} , f_{1}] = 0, \\

[e_{1} , f_{1} ]=[2 h_{1}] &
\end{array}
\end{equation}
and
\begin{equation}
\begin{array}{ll}
[ e_{2},f_{2} ]=\;[2 h_{2}]_{2}, & \\

e_{2} e_{3}^{(\pm)} = q^{\mp 2} e_{3}^{(\pm)} e_{2}, &
f_{3}^{(\pm)}  f_{2} = q^{\mp 2} f_{2} f_{3}^{(\pm)}  \\

[e_{1} , e_{4}]=0, & [f_{1} ,f_{4}]=0
\end{array}
\end{equation}
where
\begin{equation}
\begin{array}{l}
e_{3}^{(\pm)} = q^{\pm 1} e_{1} e_{2}- q^{\mp 1} e_{2} e_{1},\;\;\;\;\;\;\;
f_{3}^{(\pm)} = q^{\pm 1} f_{2} f_{1}- q^{\mp 1} f_{1} f_{2}, \\
e_{4}= q^{-1} e_{1} e_{3}^{(+)} - q\;e_{3}^{(+)} e_{1}=
q\;e_{1} e_{3}^{(-)} - q^{-1} e_{3}^{(-)} e_{1},  \\
f_{4}= q^{-1} f_{3}^{(+)} f_{1} - q\;f_{1} f_{3}^{(+)}=
q\;f_{3}^{(-)} f_{1} - q^{-1} f_{1} f_{3}^{(-)}.
\end{array}
\end{equation}
We define $(M, \;K,\;M_{2},\;M_{4})$ through
\begin{equation}
\begin{array}{ll}
q^{\pm M}=q^{\pm h_{1}}, & q^{\pm (K-M)}=q^{\pm h_{2}} \\
q^{\pm M_{2}}=q^{\pm h_{2}}, & q^{\pm M_{4}}=q^{\pm (h_{1}+h_{2})} \\
\end{array}
\end{equation}
The fundamental Casimir (classically quadratic in the Cartan-Weyl generators)
can now be written $[6,\;7]$ for arbitrary $q$ as
\begin{equation}
\begin{array}{ll}
A &= {1\over [2]} \Bigl\lbrace \bigl( f_{1}e_{1}+
[M][M+1]\bigl) {[2 K+3]_{2}\over [2 K+3]} + [K][K+3]\Bigr\rbrace \cr \\
& +\bigl(f_{2}e_{2}+{1\over [2]^{2}} f_{4}e_{4}\bigl)+
{1 \over [2]^{2}}\bigl(f_{3}^{(+)}e_{3}^{(+)}
q^{2M+1} + f_{3}^{(-)}e_{3}^{(-)}q^{-2M-1}\bigl).
\end{array}
\end{equation}

For brevity we consider in this talk only generic $q$, real positive. The case
$q$ a root of unity has been discussed in $[6]$. For generic $q$ the irreps.
are labelled by two (half) integers
\begin{equation}
\begin{array}{l}
n_{1} \geq n_{2}
\end{array}
\end{equation}
One has the following general result $[7]$. The state
annihilated by $e_{1}$, $e_{2}$ corresponds to eigenvalue $n_{2}$ of $M$ and
$n_{1}$ of $K$. Hence on this space
\begin{equation}
\begin{array}{l}
A = {1\over [2]} \Bigl\lbrace [n_{1}][n_{1}+3]+
[n_{2}][n_{2}+1]{[2n_{1}+3]_{2}\over [2n_{1}+3]}\Bigl\rbrace\;
\hbox{{\bf 1}}
\end{array}
\end{equation}
where \hbox{{\bf 1}} is the unit matrix corresponding to the dimension
\begin{equation}
\begin{array}{l}
{1\over 6}(2n_{2}+1)(2n_{1}+3)(n_{1}+n_{2}+2)(n_{1}-n_{2}+1)
\end{array}
\end{equation}
For $n_{2}=0,\;{1\over 2}$ and $n_{1}$ $(11)$ reduces to the respective
results in $[6]$. The factor $[2\;n_{1}+3]_{2} /[2\;n_{1}+3]$ in $(11)$ is a
particularly interesting example of the non-minimal case $(3)$. Its
implications will be analysed at the end.

After all the $SU(N)$ and $SO(3)$, $SO(4)$ one encounters unequal roots for
the first time for $SO(5)$. In constructing matrix elements one can associate
the well-known $SU(2)$ structure either with the Chevalley triplet
($e_{1}$, $f_{1}$, $h_{1}$) or with the triplet ($e_{2}$, $f_{2}$, $h_{2}$).
The consequences are very different and even more so concerning
$q$-deformations. Non-minimality and non-simple lacing (unequal roots) will
appear, at least in this example, associated together.

The well-known Gelfand-Zetlin basis and matrix elements $[1]$ are quite
unsuitable (for the orthogonal case) as starting point for $q$-deformation.
The situation is entirely different from that of the unitary case. The reasons
were discussed in $[6]$. After this remark I now introduce the two bases
starting with the Chevalley triplets 1 and 2 respectively.
\vskip0.5cm
{\bf Basis I.} Standard ${\cal U}_{\displaystyle{q}}(SU(2))$ structure for
($e_{1}$, $f_{1}$, $q^{\pm h_{1}}$), invariants ($n_{1}$, $n_{2}$),
variable indices ($j$, $m$, $k$, $l$). The domain of the indices are $[7]$:
\vskip 0.25cm
(i) For ($n_{1},\;n_{2}$) integers
\begin{equation}
\begin{array}{l}
j = 0,\;1,\cdots,\;n_{1}-1,\;n_{1},\;\;\;\;\;\;\;\;
m = -j,\;-j+1, \cdots,\;j-1,\;j \\
k = -l,\;-l+2, \cdots,\;l-2,\;l,\;\;\;\;\;\;\;\;
l = 0,\;1,\;2\;\cdots \\
j+l = n_{1}-n_{2},\;n_{1}-n_{2}+1, \cdots,\;n_{1}+n_{2} \\
j-l-{1\over2}\bigl(1-(-1)^{n_{1}+n_{2}-j-l}\bigl)=-n_{1}+n_{2},
\;-n_{1}+n_{2}+2, \cdots,\;n_{1}-n_{2}.
\end{array}
\end{equation}
\vskip 0.25cm
(ii) For ($n_{1},\;n_{2}$) half-integers
\begin{equation}
\begin{array}{l}
j = {1\over 2},\;{3\over 2},\cdots,\;n_{1}-1,\;n_{1},\;\;\;\;\;\;\;\;
m = -j,\;-j+1, \cdots,\;j-1,\;j \\
k = -l,\;-l+1, \cdots,\;l-1,\;l,\;\;\;\;\;\;\;\;
l = {1\over 2},\;{3\over 2},\cdots \\
j+l = n_{1}-n_{2}+1,\;n_{1}-n_{2}+3, \cdots,\;n_{1}+n_{2} \\
j-l=-n_{1}+n_{2},\;-n_{1}+n_{2}+2, \cdots,\;n_{1}-n_{2}.
\end{array}
\end{equation}
The dimension is given by $(12)$ for both cases.

The matrix elements are (suppressing $n_{1}$, $n_{2}$ in the state labels)
\begin{equation}
\begin{array}{l}
q^{\pm M } \vert j\;m\;k\;l\rangle
= q^{\pm m } \vert j\;m\;k\;l\rangle \\
q^{\pm K } \vert j\;m\;k\;l\rangle
= q^{\pm k} \vert j\;m\;k\;l\rangle \\
e_{1} \vert j\;m\;k\;l\rangle
= ([j-m]\;[j+m+1])^{1/2}\vert j\;m+1\;k\;l \rangle \\
e_{2} \vert j\;m\;k\;l\rangle=
([j-m+1][j-m+2])^{1/2}\sum_{ l'}\;a(j,k,l,l')
\vert j+1\;m-1\;k+1\;l'\rangle \\
\;\;\;\;\;\;\;\;\;\;\;\;\;\;\;\;\;\;+ (\;[j+m][j+m-1])^{1/2}\sum_{
l'}\;b(j,k,l,l')
\vert j-1\;m-1\;k+1\;l'\rangle \\
\;\;\;\;\;\;\;\;\;\;\;\;\;\;\;\;\;\;+(\;[j+m][j-m+1])^{1/2}\sum_{
l'}\;c(j,k,l,l')
\vert j\;m-1\;k+1\;l'\rangle \\
\end{array}
\end{equation}
We consider only real solutions of the matrix elements when for any two
states $|x\rangle$, $|y\rangle$
\begin{equation}
\begin{array}{l}
\langle y | f_{i} |x\rangle = \langle x | e_{i} |y\rangle,\;\;\;\;\;\;\;
(i=1,\;2)
\end{array}
\end{equation}
As yet solutions have been obtained $[6]$ for arbitrary (half) integer $n_{1}$
only for the extreme values of $n_{2}$,
\begin{equation}
\begin{array}{l}
n_{2}=0\;\;\;\hbox{or}\;\;\;{1\over 2}
\end{array}
\end{equation}
and
\begin{equation}
\begin{array}{l}
n_{2}=n_{1}
\end{array}
\end{equation}
For these cases $l$-dependence is trivial. One labels the states as
$|j\;m\;k\rangle$. Even classical $(q=1)$ solutions are not available for
the general case. Referring to $[6]$ for complete solutions of the cases
$(17)$ I now present the solution of $(18)$ for comparing it to the
corresponding one in Basis II to follow.

For $n_{2}=n_{1}=n$ (integer or half-integer), suppressing trivial
$l$-dependence
\begin{equation}
\begin{array}{l}
a(j,k) =b(j+1,-k-1)=(q+q^{-1})^{-1} \Biggl({[n-j]_{2}\;[n+j+2]_{2}\;[j+k+1]\;
[j+k+2]\over [2j+3]\;[2j+1]\;[j+1]_{2}^{2}}\Biggl)^{1/2} \\
c(j,k) =(q+q^{-1})^{-1} [n+1]_{2} {([j-k]\;[j+k+1])^{1/2} \over [j+1]_{2}
\;[j]_{2}}
\end{array}
\end{equation}
The dimension is
\begin{equation}
\begin{array}{l}
{1\over 3}(n+1)(2n+1)(2n+3)
\end{array}
\end{equation}

Comparing with the limit $q=1$, it is evident that each factor undergoes a
(pseudo) minimal deformation of type $(1)$. The situation will change in the
basis to follow.
\vskip 0.5cm
{\bf Basis II.} Standard ${\cal U}_{q^{2}}(SU(2))$ structure
for ($e_{2}$, $f_{2}$, $q^{\pm h_{2}}$), invariants ($n_{1}$, $n_{2}$),
variable indices ($j_{2}$, $m_{2}$, $j_{4}$, $m_{4}$). The domain of the
indices are $[7]$ (for integer and half-integer $n_{1}$, $n_{2}$):
\begin{equation}
\begin{array}{l}
j_{2} = 0,\;{1\over 2},\;1,\cdots,\;{n_{1}+n_{2}\over 2},\;\;\;\;\;\;\;\;
m_{2} = -j_{2},\;-j_{2}+1, \cdots,\;j_{2}-1,\;j_{2} \\
j_{4} = 0,\;{1\over 2},\;1,\cdots,\;{n_{1}+n_{2}\over 2},\;\;\;\;\;\;\;\;
m_{4} = -j_{4},\;-j_{4}+1, \cdots,\;j_{4}-1,\;j_{4} \\
j_{2}+j_{4} = n_{2},\;n_{2}+1,\cdots,\;n_{1} \\
j_{2}-j_{4} = -n_{2},\;-n_{2}+1,\cdots,\;n_{2}.
\end{array}
\end{equation}
The dimension is given by $(12)$.

The matrix elements are
\begin{equation}
\begin{array}{l}
q^{\pm M_{2}} \vert j_{2}\;m_{2}\;j_{4}\;m_{4}\rangle
= q^{\pm m_{2}} \vert j_{2}\;m_{2}\;j_{4}\;m_{4}\rangle \\
q^{\pm M_{4}} \vert j_{2}\;m_{2}\;j_{4}\;m_{4}\rangle
= q^{\pm m_{4}} \vert j_{2}\;m_{2}\;j_{4}\;m_{4}\rangle \\
e_{2} \vert j_{2}\;m_{2}\;j_{4}\;m_{4} \rangle
= ([j_{2} - m_{2}]_{2} [j_{2}+m_{2}+1]_{2})^{1/2}
\vert j_{2}\;m_{2}+1\;j_{4}\;m_{4} \rangle \\
e_{1} \vert j_{2}\;m_{2}\;j_{4}\;m_{4}\rangle=
\sum_{\epsilon,\; \epsilon '}\;\;
([j_{2}-\epsilon\;m_{2}+{1+\epsilon \over 2}]_{2})^{1/2}\; \\
\;\;\;\;\;\;\;\;\;\;\;\;\;\;\;\;\;\;\;\;\;\;\;\;\;\times\;
c_{(\epsilon,\epsilon ')}(j_{2},j_{4},m_{4})
\;\vert j_{2}+{\epsilon \over 2}\;\;m_{2}-{1\over 2}\;\;j_{4}+
{\epsilon ' \over 2}\;\;m_{4}+{1\over 2}\rangle
\end{array}
\end{equation}
with, as in $(16)$,
\begin{equation}
\begin{array}{l}
\langle y | f_{i} |x\rangle = \langle x | e_{i} |y\rangle,\;\;\;\;\;\;\;
(i=1,\;2)
\end{array}
\end{equation}

Now a classical solution is available $[8]$. In our notation this is
\begin{equation}
\begin{array}{l}
c_{(\epsilon,\epsilon ')}(j_{2},j_{4},m_{4})=(j_{4}+\epsilon '\;m_{4}+
{1+\epsilon ' \over 2})^{1/2} c_{(\epsilon,\epsilon ')}(j_{2},j_{4})\;\;\;
(\epsilon,\epsilon '=\pm 1) \\
c_{(\epsilon,\epsilon ')}(j_{2},j_{4})=\epsilon\; \epsilon '\;
c_{(-\epsilon,-\epsilon ')}(j_{2}+{\epsilon \over 2},j_{4}+
{\epsilon ' \over 2})
\end{array}
\end{equation}
where
\begin{equation}
\begin{array}{l}
c_{(++)}(j_{2},\;j_{4})=\biggl({(n_{1}+j_{2}+j_{4}+3)(n_{1}-j_{2}-j_{4})
(j_{2}+j_{4}+n_{2}+2)(j_{2}+j_{4}-n_{2}+1) \over (2j_{2}+1)\;(2j_{2}+2)\;
(2j_{4}+1)\;(2j_{4}+2)}\biggl )^{1/2} \\
c_{(+-)}(j_{2},\;j_{4})=\biggl({(n_{1}+j_{2}-j_{4}+2)(n_{1}-j_{2}+j_{4}+1)
(j_{2}-j_{4}+n_{2}+1)(j_{4}-j_{2}+n_{2}) \over (2j_{2}+1)\;(2j_{2}+2)\;
(2j_{4})\;(2j_{4}+1)}\biggl )^{1/2}.
\end{array}
\end{equation}
(see the comments in $[7]$ concerning the relation to $[8]$).

But now $q$-deformation is the problem. So far solutions have been obtained for
\begin{equation}
\begin{array}{l}
n_{2}=0
\end{array}
\end{equation}
and
\begin{equation}
\begin{array}{l}
n_{2}=n_{1}=n
\end{array}
\end{equation}
Referring to $[7]$ for $(26)$ I now reproduce only the solution for $(27)$.

For $n_{1}=n_{2}=n$
\begin{eqnarray*}
j_{2}+j_{4}=n,\;\;\;\;\;\;\;\;\;\;\;\;\;\;\;&j_{2}=0,\;{1\over 2},\cdots,\;n \\
c_{(\epsilon \epsilon)}(j_{2},\;j_{4},\;m_{4})=0.  &  \\
\end{eqnarray*}
and
\begin{eqnarray*}
c_{(\epsilon , -\epsilon)}(j_{2},\;j_{4},\;m_{4})=([n+1]_{2}-
[j_{2}+\epsilon\;m_{4}+{1 \over 2}(1+\epsilon)]_{2})^{1/2}\;
c_{(\epsilon , -\epsilon)}(j_{2})
\end{eqnarray*}
with
\begin{equation}
\begin{array}{l}
c_{(+-)}(j_{2})=-c_{(-+)}(j_{2}+{1\over 2})=\biggl({[2j_{2}+1]\;[2j_{2}+2]\over
[2j_{2}+1]_{2}\;[2j_{2}+2]_{2}}\biggl)^{1/2}.
\end{array}
\end{equation}
The $m_{4}$-dependence in $(28)$ is of type $(2)$ with
$$
x_{1}=n+1,\;\;\;\;x_{2}=j_{2}+\epsilon m_{4}+{1 \over 2}(1+\epsilon)
$$
so that
\begin{equation}
\begin{array}{l}
x=x_{1}-x_{2}=j_{4}+\epsilon ' m_{4} +{1\over 2}(1+\epsilon')
\end{array}
\end{equation}
consistensly with $(24)$. The $c(j_{2})$'s in $(28)$ are of type $(3)$ with
$x=1$ but double $y$-factors.

{\it If one tries to solve using only (pseudo) minimal deformations of type
$(1)$ one runs into contradictions}. Thus non-minimality is essential for this
basis. This basis, in turn, seems to be essential for providing access to
certain interesting sectors. Let me just mention two such points.

(i) Suitably adapting familiar continuation techniques Basis I leads to
${\cal U}_{q}(SO(3,2))$ while Basis II is needed to arrive at
${\cal U}_{q}(SO(4,1))$.

(ii) Under suitable contraction procedures quite different $q$-deformed
inhomogeneous algebras are obtained from the two bases respectively (see the
comments in $[6]$ and $[7]$).

It is not possible to discuss these aspects here. But they suffice to indicate
the potential interest of a general solution of $(22)$ for arbitrary $q$. (In
fact once solutions are found for generic $q$ our method of fractional parts
explained in $[3]$ and $[4]$ and already used for Basis I in $[6]$ will
readily yield solutions for $q$ a root of unity.)

The general solution will also permit a better understanding of the role of
non-minimality. This role is not merely formal. If there is indeed some
physical application, the physical content of different types of deformations
will be different. As $q$ moves away from unity they will respond differently.
Thus, to take an example, for $q=e^{\delta}$ and $x=x_{1}-x_{2}$
\begin{equation}
\begin{array}{l}
[x_{1}]_{p_{1}}-[x_{2}]_{p_{2}}=x+{1\over 6}\delta^{2}(p_{1}^{2}
x_{1}(x_{1}^{2}-1)-p_{2}^{2}x_{2}(x_{2}^{2}-1))+\cdots
\end{array}
\end{equation}
In this context the non-minimality of type $(3)$ noted in $(11)$ also has a
striking consequence.

For $q=e^{\delta}$,
\begin{equation}
\begin{array}{l}
{1\over [2]} \Bigl\lbrace [n_{1}][n_{1}+3]+
[n_{2}][n_{2}+1]{[2n_{1}+3]_{2}\over [2n_{1}+3]}\Bigl\rbrace=A_{2}+
\delta^{2}\;A_{4}+\cdots
\end{array}
\end{equation}
where
\begin{equation}
\begin{array}{l}
A_{2}={1\over 2}(n_{1}(n_{1}+3)+n_{2}(n_{2}+1))
\end{array}
\end{equation}
\begin{equation}
\begin{array}{l}
A_{4}=4\;n_{2}(n_{2}+1)(n_{1}+1)(n_{1}+2)+\cdots
\end{array}
\end{equation}
The other terms of $A_{4}$ are very easily obtained. Let us, however,
concentrate on the first term, a direct consequence of the factor
$[2n_{1}+3]_{2}/ [2n_{1}+3]$ in $(31)$.

$A_{2}$ is just the well-known eigenvalue of the first classical Casimir
(quadratic in the Cartan-Weyl generators) for the irrep. ($n_{1}$, $n_{2}$).
The first term of $A_{4}$ is the classical eigenvalue of the second
(quadratic) Casimir operator.

{\it This is an illustration of the general result announced in my first talk
$[9]$. The $q$-deformed quadratic Casimir alone completely characterizes the
irreps. ($n_{1}$, $n_{2}$). We note morever that this is here achieved through
a typical non-minimality.}

Consider now the consequence of the same non-minimal structure in the context
of contraction. Certain aspects of contraction of ${\cal U}_{q}(SO(5))$ are
discussed in $[6]$. Here let us just note that the eigenvalue of the
contracted Casimir ( for $q > 1$ for example ) is obtained by dividing the
l.h.s. of $(31)$ by the leading term of $[n_{1}][n_{1}+3]$ as
$n_{1}\rightarrow \infty$ multiplied by a constant $\lambda^{-2}$, i.e. by
\begin{equation}
\begin{array}{l}
{1 \over \lambda^{2}[2]}{q^{2n_{1}+3} \over (q-q^{-1})^{2}}
\end{array}
\end{equation}
and taking the limit. This gives an eigenvalue
\begin{equation}
\begin{array}{l}
\lambda^{2}\lbrace 1+ {(q-q^{-1})^{2} \over (q+q^{-1})}[n_{2}][n_{2}+1]
\rbrace
\end{array}
\end{equation}

A general solution for ${\cal U}_{q}(SO(5))$ on a suitable basis can lead
through contraction to a successful construction of representations of
${\cal U}_{q}(E(4))$ (the $q$-deformed $4$-dimensional Euclidean algebra) for
arbitrary $q$. Then one has to see if a suitable analytic continuation to
$q$-deformed Poincar\'e algebra is possible through this approach. This is one
of our main goals. This remains to be done. {\it But (35) shows that it will
include the following remarkable feature. The $q$-deformed "mass-like"
operator (the $q$-deformation of the sum of squares of the translations) will
have eigenvalues depending on the "spin-like" parameter $n_{2}$ as in $(35)$}.
This is reminiscent of a famous feature of $SU(6)$ type models. This seeping
of internal discrete quantum numbers into the "mass-like" spectrum seems to
be a typical feature of $q$-deformations $[10]$. But for the orthogonal
symmetry (at least for the present example) this turns out to be an intriguing
consequence of non-minimality. One need not inject ans${\ddot{a}}$tze to
construct a mass spectrum depending on internal quantum numbers. One just
solves the mathematical problem of constructing representations explicitly
and the result is there.

\end{document}